\documentclass[11pt]{article}   
\usepackage{amsthm, amsmath}   
\usepackage{amssymb,amsfonts,amsbsy}
\usepackage[latin1]{inputenc}

\newcommand{\Q}{\mathbb{Q}}
\newcommand{\C}{\mathbb{C}}
\newcommand{\Z}{\mathbb{Z}}
\newcommand{\R}{\mathbb{R}}
\newcommand{\N}{\mathbb{N}}
\newcommand{\M}{\mathbb{M}}
\newcommand{\K}{\mathbb{K}}
\newcommand{\gr}{{\rm gr}}

\theoremstyle{plain}
\newtheorem{teorema}{Teorema}[section]
\newtheorem{lema}[teorema]{Lema}

\theoremstyle{definition}
\newtheorem{definicion}[teorema]{Definici\'on}

\newtheorem{algoritmo}[teorema]{Algoritmo}

\theoremstyle{remark}

\title{Aplicación de la descomposición racional univariada a Monstrous Moonshine}
\author{John McKay \and David Sevilla}
\date{}

\begin{document}
\maketitle

\begin{abstract}
This paper shows how to use Computational Algebra techniques,
namely the decomposition of rational functions in one variable, to
explore a certain set of modular functions, called replicable
functions, that arise in Monstrous Moonshine. This is an active
topic in Mathematics that links the Monster group and complex
analysis due to a surprising relation between the character table
of the Monster and the coefficients of the expansions of certain
modular functions. In particular, we have computed all the
rational relations with coefficients in $\Z$ between pairs of
replicable functions.
\end{abstract}

\section{Introducci\'on}

Las funciones modulares son objetos clásicos en teoría de números,
y resultan de considerar la acción del grupo $PSL_2(\Z)$ sobre el
plano hiperbólico (representado, por ejemplo, como el semiplano
complejo superior). El estudio de ciertos grupos de
transformaciones, conmensurables con $PSL_2(\Z)$ (es decir,
aquellos cuya intersección tiene índice finito en ambos)
proporciona una clase de funciones que quedan invariantes por
ellos, de las que la clásica función $j$ es representante. En
concreto, los normalizadores de todos estos grupos en $PSL_2(\R)$
del subgrupo de matrices triangulares superiores módulo $N$.

En los años 70, McKay hizo notar que los coeficientes de ciertos
desarrollos en serie de estas funciones (que son de hecho enteros
positivos) est\"an relacionados con las entradas de la tabla de
caracteres del grupo monstruo $\M$. Conway y Norton (\cite{CN79})
formularon diversas conjeturas sobre estas funciones, que fueron
probadas entre otras por Borcherds en \cite{Bor92}. A continuación
presentamos brevemente algunas definiciones y resultados
conocidos.

\begin{definicion}
Dados $\omega_1,\omega_2$ períodos de una función doblemente
periódica con $\tau=\omega_2/\omega_1\in\mathcal{H}=\{z\in\C:\
Im(z)>0\}$, la función $J$ (conocida comúnmente como función
modular de Klein) es
\[J(\omega_1,\omega_2)=\frac{g_2^3(\omega_1,\omega_2)}{\Delta(\omega_1,\omega_2)},\]
donde $g_2,g_3$ son los invariantes de la función elíptica de
Weierstrass con discriminante $\Delta=g_2^3-27g_3^2$.

Para cada $\tau\in \mathcal{H}$, definimos
$J(\tau)=J(1,\tau)=J(\omega_1,\omega_2)$.

La función $j$ se define clásicamente como
\[j(q)=1728\cdot J(\sqrt{q}).\]

Los primeros coeficientes de $j$ son:
\[j(q)=\frac{1}{q}+744+196884\,q+21493760\,q^2+864299970\,q^3+20245856256\,q^4+\cdots\]
\end{definicion}

\begin{definicion}
El \emph{grupo modular} es el grupo de transformaciones
\[z\rightarrow\frac{az+b}{cz+d}\]
donde $a,b,c,d\in\Z, ad-bc=1$, representado proyectivamente por
las matrices $2x2$ con coeficientes en $\Z$.
\end{definicion}

Las funciones con las que trabajamos son aquellas fijadas por
grupos discretos $G_f$, es decir,
\[f(\gamma z)=f(z) \forall\gamma\in G_f.\]

En general, hay relaciones polinomiales entre cada dos de estas
funciones. En particular, cuando uno de los grupos está contenido
en el otro se tiene que una de las funciones se puede escribir
como función racional de la otra (si los grupos son simplemente
conmensurables, puede ser necesario tomar una potencia de $q$ en
la primera función).

En particular, las funciones replicables surgen como
generalización de las funciones modulares relacionadas con el
grupo $\M$ (ver \cite{ACMS92}). Nuestra aportación consiste en
intentar refinar al máximo el poset de estas 619 funciones,
encontrando relaciones entre ellas del tipo $s_1=f(s_2)$ con
$f\in\C(x)$ o $f\in\Q(x)$, donde $s_1, s_2$ son series en
$q^{k_1}, q^{k_2}$, y descomponiendo estas funciones racionales
para encontrar cadenas maximales en ese poset.

Nuestros cálculos permiten, entre otras cosas, comprobar si la
lista de 619 funciones contiene todas las funciones replicables, o
si aparecen funciones de otras características especiales como
resultado de este refinamiento.

\section{Algoritmos}

Tenemos el siguiente algoritmo elemental para encontrar relaciones
racionales entre dos series en general.

\begin{algoritmo} $ $
\begin{description}
    \item \textsc{Entrada}: dos funciones replicables $j_1$ y $j_2$.
    \item \textsc{Salida}: una función $f\in\Q(x)$ tal que
$j_1(q^r)=f(j_2(q))$ para algún $r\in\N$, si existe.
\end{description}
\begin{description}
    \item \textsc{A}. Calcular las áreas $A_1$, $A_2$ de las regiones
fundamentales de $j_1$ y $j_2$ respectivamente. Si $e=A_2/A_1$ no
es un número natural, terminar. Si lo es, sea $r=1$.

    \item \textsc{B}. Sea
\[f=\frac{t^e+a_{e-1}t^{e-1}+\cdots+a_0}{t^{e-r}+b_{e-r-1}t^{e-r-1}+\cdots+b_0}.\]
Sean $s_1=1/q+\sum_{k=0}^{2e+1}c_kq^k$,
$s_2=1/q+\sum_{k=0}^{2e+1}d_kq^k$ las series truncadas de $j_1$ y
$j_2$ respectivamente.

    \item \textsc{C}. Resolver el sistema de ecuaciones lineales
en las variables $a_i,b_j$ dado por la anulación del numerador de
$j_1(q^r)-f(j_2(q))$. Si tiene solución, devolver el $f$
correspondiente.

    \item \textsc{D}. Si $r<e$, incrementar $r$ y volver a
\textsc{B}. Si no, terminar.
\end{description}
\end{algoritmo}

Este es un algoritmo de complejidad baja, ya que las áreas de las
regiones fundamentales son conocidas (fácilmente computables a
partir de generadores de los grupos que las fijan) y calcular los
coeficientes se hace a través de la resolución de sistemas de
ecuaciones lineales. De esta manera construimos un grafo cuyos
nodos son funciones modulares y cuyas aristas vienen dadas por
pares de funciones relacionadas de la manera anterior.

El siguiente paso es utilizar un algoritmo propio de
descomposición de funciones racionales univariadas, eficiente en
la práctica, para refinar todo lo posible este grafo. Este
algoritmo, cuyos detalles pueden encontrarse en \cite{Sev04}, está
basado en una idea presentada en \cite{AGR95}. Damos a
continuación las definiciones y resultados más relevantes.

\begin{definicion}
Definimos el \emph{grado} de una función racional como el máximo
de los grados de su numerador y su denominador, suponiendo que no
tienen factores comunes.

Una función racional $f\in\K(x)$ está en \emph{forma normal} si
$\gr\ f_N>\gr\ f_D$ y $f_N(0)=0$ (por tanto $f_D(0)\neq 0$).
\end{definicion}

\begin{lema}\label{exist-unids}
Dada $f\in\K(x)$, si $\gr\ f<|\K|$ existen unidades $u,v$ tales
que $u\circ f\circ v$ está en forma normal.
\end{lema}

\begin{teorema}\label{desc-equiv-normal}
Si $f$ está en forma normal, toda descomposición suya es
equivalente a una descomposición en la que ambas componentes están
en forma normal.
\end{teorema}

\begin{teorema}\label{divide-normal-univ}
Sea $f=g(h)$ con $f,g,h$ en forma normal. Entonces $h_N$ divide a
$f_N$ y $h_D$ divide a $f_D$.
\end{teorema}

El teorema anterior nos proporciona el algoritmo de descomposición
antes mencionado.

\begin{algoritmo}\label{alg-desc-univ} $ $
\begin{description}
    \item \textsc{Entrada}: $f\in\K(x)$.
    \item \textsc{Salida}: todas las descomposiciones de $f$ (es decir, al menos
    un representante de cada clase de equivalencia de descomposiciones).
\end{description}
\begin{description}
    \item \textsc{A}. Calcular unidades $u,v$ como en el Lema
    \ref{exist-unids}. Sea $\overline{f}=u\circ f\circ v$.

    \item \textsc{B}. Factorizar $\overline{f}_N$ y $\overline{f}_D$. Sea
    $D=\{(A_1,B_1),\ldots,(A_m,B_m)\}$ el conjunto de pares $(A,B)$ de polinomios
    mónicos tales que $A,B$ dividen a $\overline{f}_N,\overline{f}_D$ respectivamente.

    \item \textsc{C}. Para cada $i\in\{1,\ldots,m\}$ comprobar si existe $g\in\K(x)$ tal que
    $\overline{f}=g\circ (A_i/B_i)$. Si es así, añadir
    $\left(u^{-1}(g),h(v^{-1})\right)$ a la lista de
    descomposiciones.

    \item \textsc{D}. Devolver la lista de descomposiciones (si no se ha
    encontrado ninguna, $f$ es indescomponible).
\end{description}
\end{algoritmo}

\section{Resultados}

La utilización de esta combinación de técnicas nos ha permitido
calcular este grafo refinado descomponiendo las relaciones en
$\Q(x)$ encontradas de acuerdo con este esquema:

\[j_1\ \stackrel{d,r}{\longrightarrow}\ j_2\ \ \ \Rightarrow\ \ \ j_1(q^r)=f(j_2(q)),\ \gr\:f=d\]

\[f=g\circ h\ \ \ \Rightarrow\ \ \ \left\{\begin{array}{l}
j_3(q^s)=h(j_2(q)) \\
\\
j_1(q^r)=g(j_3(q^s)) \\
\end{array}\right.\]

\[j_1\ \stackrel{\gr\:g,r/s}{\longrightarrow}\ j_3\ \stackrel{\gr\:h,s}{\longrightarrow}\ j_2\]

En el primer paso se obtuvieron 2419 relaciones. Se necesitaron
aproximadamente 20 horas en un ordenador personal PC Pentium 4
utilizando el software de cálculo simbólico Maple para completar
esta etapa; la mayor parte del tiempo se utilizó en el cálculo de
coeficientes de las series involucradas por medio de una
recurrencia. En la siguiente tabla mostramos los grados de las
funciones racionales encontradas:

\[\begin{array}{ccccc}
grado & cantidad & \ \ \ & grado & cantidad \\
2 & 698 & & 16 & 52 \\
3 & 243 & & 18 & 60 \\
4 & 422 & & 20 & 2 \\
5 & 26 & & 24 & 71 \\
6 & 333 & & 28 & 2 \\
8 & 178 & & 30 & 8 \\
9 & 40 & & 32 & 4 \\
10 & 14 & & 36 & 40 \\
12 & 209 & & 48 & 5 \\
14 & 4 & & 72 & 2 \\
15 & 6 \\
\end{array}\]

En la etapa de descomposición, tras 30 horas de cálculo se computó
el grafo buscado, que consta de 619 nodos y 1202 aristas.

El interés de los resultados obtenidos es múltiple. Por una parte,
las relaciones obtenidas son una buena comprobación de que las
tablas de datos conocidas son correctas, dado que no se han
encontrado funciones nuevas mediante este método.

Por otra parte, también se han observado interesantes propiedades
en algunas de las funciones racionales encontradas. El ejemplo más
claro en ese sentido es la relación racional entre las funciones
denotadas $(1A)$ y $(9B)$:
\[f=\frac{x^3(x+6)^3(x^2-6\,x+36)^3}{(x-3)^3(x^2+3\,x+9)^3}.\]

Durante el último paso fueron encontradas dos descomposiciones
completas de $f$:
\[f=x^3\circ\frac{x(x-12)}{x-3}\circ\frac{x(x+6)}{x-3}=\frac{x^3(x+24)}{x-3}\circ\frac{x(x^2-6\,x+36)}{x^2+3\,x+9}.\]

Esto significa que la función $f$ tiene dos cadenas completas de
descomposición de distinta longitud. Por lo que sabemos, es el
primer ejemplo conocido de función racional en $\Q(x)$ con esta
propiedad; como contraste, es conocido que todas las cadenas
completas de descomposición de polinomios con coeficientes en
característica cero tienen la misma longitud (ver \cite{Rit22},
\cite{Sch00}). Para más detalles al respecto, ver \cite{Sev04}.

También es interesante notar que, en algunos casos, para dos
series $s_1,s_2$ encontramos más de una relación, es decir,
$s_1(q)=f_1(s_2(q^{k_1}))=f_2(s_2(q^{k_2}))$ con $k_1\neq k_2$.
Con una simple resultante podemos encontrar una relación
polinomial del tipo $P(s_2(q),s_2(q^k))=0$. Estas relaciones
tienen un interés propio en la teoría relativa al Monstrous
Moonshine, siendo generalizaciones de los clásicos
\emph{polinomios modulares} para $j$.

Una posible línea de investigación inmediata consiste en calcular
todas las descomposiciones de estas funciones en $\C(x)$, dado que
se espera que las acciones no triviales que surgen sobre los
grupos correspondientes proporcionen nuevas funciones de interés.
También esperamos poder aplicar otras técnicas de Álgebra
Computacional a aspectos relacionados dentro del Monstrous
Moonshine.

\vspace{2em}

David Sevilla González

Dpto. de Matemáticas, Estadística y Computación - Facultad de
Ciencias

Universidad de Cantabria

Av. de los Castros s/n

39071 Santander (Cantabria)

david.sevilla@unican.es

\vspace{2em}

John McKay

Sir George Williams Campus

1455 De Maisonneuve Blvd. West

LB 903 15

Montreal, QC, Canada, H3G 1M8

mckay@cs.concordia.ca

\end{document}